\documentclass[12pt]{article}
\usepackage{amssymb}
\usepackage{latexsym,bm}
\usepackage{graphicx}
\usepackage{amsmath}
\usepackage{mathrsfs}

\newcommand{\bea}{\begin{eqnarray*}}
\newcommand{\eea}{\end{eqnarray*}}
\newcommand{\be}{\begin{equation}}
\newcommand{\ee}{\end{equation}}
\newcommand{\ben}{\begin{eqnarray*}}
\newcommand{\een}{\end{eqnarray*}}

\date{}
\begin{document}
\title{On a problem of Erd\H{o}s about graphs whose size is the Tur\'{a}n number plus one\footnote{E-mail addresses:
{\tt pq@ecust.edu.cn}(P.Qiao),
{\tt zhan@math.ecnu.edu.cn}(X.Zhan).}}
\author{\hskip -10mm Pu Qiao$^a$, Xingzhi Zhan$^b$\thanks{Corresponding author.}\\
{\hskip -10mm $^a$\small  Department of Mathematics, East China University of Science and Technology, Shanghai 200237, China}\\
{\hskip -10mm $^b$\small Department of Mathematics, East China Normal University, Shanghai 200241, China}}\maketitle
\begin{abstract}
 We consider finite simple graphs. Given a graph $H$ and a positive integer $n,$ the Tur\'{a}n number of $H$ for the order $n,$ denoted
 ${\rm ex}(n,H),$ is the maximum size of a graph of order $n$ not containing $H$ as a subgraph. Erd\H{o}s posed the following problem in 1990:
 ``For which graphs $H$ is it true that every graph on $n$ vertices and ${\rm ex}(n,H)+1$ edges contains at least two $H$s? Perhaps this is always true."
 We solve the second part of this problem in the negative by proving that for every integer $k\ge 4,$ there exists a graph $H$ of order $k$ and
 at least two orders $n$ such that there exists a graph of order $n$ and size ${\rm ex}(n,H)+1$ which contains exactly one copy of $H.$
 Denote by $C_4$ the $4$-cycle. We also prove that for every integer $n$ with $6\le n\le 11,$ there exists a graph of order $n$ and
 size ${\rm ex}(n,C_4)+1$ which contains exactly one copy of $C_4,$ but for $n=12$ or $n=13,$ the minimum number of copies of $C_4$
 in a graph of order $n$ and size ${\rm ex}(n,C_4)+1$ is $2.$ 
\end{abstract}

{\bf Key words.} Tur\'{a}n number; extremal graph theory; Erd\H{o}s' problem; book; star; 4-cycle

{\bf Mathematics Subject Classification.} 05C35, 05C30, 05C75

\section{Introduction and Statement of the Main Results}

We consider finite simple graphs, and use standard terminology and notations. The {\it order} of a graph is its number of vertices, and the
{\it size} its number of edges. Denote by $V(G)$ and $E(G)$ the vertex set and edge set of a graph $G$ respectively. For graphs we will use
equality up to isomorphism, so $G_1=G_2$ means that $G_1$ and $G_2$ are isomorphic. Given graphs $H$ and $G,$ a {\it copy} of $H$ in $G$ is
a subgraph of $G$ that is isomorphic to $H.$  $K_p$ and $C_p$ denote
the complete graph of order $p$ and cycle of length $p$ respectively. $K_{s,\, t}$ denotes the complete bipartite graph on $s$ and $t$ vertices.
In particular, $K_{1,\, p}$ is the {\it star} of order $p+1.$ A {\it triangle-free} graph is one that contains no triangles.

In 1907, Mantel [10] proved that the maximum size of a triangle-free graph of order $n$ is $\lfloor n^2/4\rfloor$ and the balanced complete
bipartite graph $K_{\lfloor n/2\rfloor,\lceil n/2\rceil}$ is the unique extremal graph. Later in 1941, Tur\'{a}n [12] solved the corresponding
problem with the triangle replaced by a general complete graph.

{\bf Definition 1.} Given a graph $H$ and a positive integer $n,$ the {\it Tur\'{a}n number of $H$ for the order $n,$} denoted ${\rm ex}(n,H),$
is the maximum size of a simple graph of order $n$ not containing $H$ as a subgraph.

Thus, Mantel's theorem says that ${\rm ex}(n,K_3)=\lfloor n^2/4\rfloor,$ and  Tur\'{a}n determined ${\rm ex}(n,K_p).$
Determining the Tur\'{a}n number for various graphs $H$ is one of the main topics in extremal graph theory [1]. Note that a triangle is both
$K_3$ and $C_3.$ It is natural to extend Mantel's theorem to the case of larger cycles. This is difficult for even cycles. For example,
the values ${\rm ex}(n,C_4)$ have not been determined for a general order $n$ (precise values are known only for some orders
of special forms; see [6] and [7]), and a conjecture of Erd\H{o}s and Simonovits on ${\rm ex}(n,C_6)$ was refuted in [8].

On the other hand, Rademacher (orally to Erd\H{o}s who gave a simple proof in 1955 [5]) proved that every graph of order $n$ and size
${\rm ex}(n,K_3)+1$ contains at least $\lfloor n/2\rfloor$ triangles. A similar result for a general $K_p$ was proved by Moon [11].

In 1990, Erd\H{o}s posed the following problem in his paper entitled ``Some of my favourite unsolved problems" [4,  pp.472-473]:

{\bf Problem 1.} For which graphs $H$ is it true that every graph on $n$ vertices and ${\rm ex}(n,H)+1$ edges contains at least two $H$s?
Perhaps this is always true.

We solve the second part of this problem in the negative. We will use the class of book graphs.

{\bf Definition 2.} The {\it book} with $p$ pages, denoted $B_p$, is the graph that consists of $p$ triangles sharing a common edge.

The main results are as follows.

{\bf Theorem 1.} {\it Let $p$ be an even integer and let $n$ be an odd integer with $n\ge p+1\ge 5.$ Then there exists a graph of order $n$
and size ${\rm ex}(n,K_{1,\, p})+1$ which contains exactly one copy of $K_{1,\, p}.$}

{\bf Theorem 2.} {\it Let $p$ be an even positive integer. Then there exists a unique graph of order $p+2$ and size ${\rm ex}(p+2,B_p)+1$
which contains exactly one copy of $B_p,$ and there exists a unique graph of order $p+3$ and size ${\rm ex}(p+3,B_p)+1$
which contains exactly one copy of $B_p.$ }

Combining Theorems 1 and 2 we obtain the following corollary.

{\bf Corollary 3.} {\it For every integer $k\ge4,$ there exists a graph $H$ of order $k$ and at least two orders $n$ such that
there exists a graph of order $n$ and size ${\rm ex}(n,H)+1$ which contains exactly one copy of $H.$ }

We remark that the conclusion in Theorem 1 is false for odd $p\ge 3$ and for the case when both $p$ and $n$ are even. In these two cases,
the minimum number of copies of $K_{1,\, p}$ in a graph of order $n$ and size ${\rm ex}(n,K_{1,\, p})+1$ is $2.$

The following result shows that the statement for odd $p$ on books corresponding to Theorem 2 is false.

{\bf Theorem 4.} {\it Let $p\ge 3$ be an odd integer. Then the minimum number of copies of $B_p$ in a graph of order $p+2$
and size ${\rm ex}(p+2,B_p)+1$  is $3,$ and the minimum number of copies of  $B_p$ in a graph of order $p+3$
and size ${\rm ex}(p+3,B_p)+1$  is $3(p+1).$ }

Recently He, Ma and Yang [9, Conjecture 10.2] made the conjecture that ${\rm ex}(q^2+q+2,C_4)=(q(q+1)^2)/2+2$ for large $q=2^k.$
They [9, Proposition 10.3] proved that if this conjecture is true, then the 4-cycle $C_4$ would serve as a counterexample to
Erd\H{o}s' Problem 1 and that [9, p.38] $C_4$ for the order $22$ is such an example. Our next result shows that $C_4$ is a 
counterexample to Problem 1 for several low orders.

{\bf Theorem 5.} {\it For every integer $n$ with $6\le n\le 11,$ there exists a graph of order $n$ and size ${\rm ex}(n,C_4)+1$
which contains exactly one copy of $C_4,$ but for $n=12$ or $n=13,$ the minimum number of copies of $C_4$ in a graph of order $n$
and size ${\rm ex}(n,C_4)+1$ is $2.$  }

In Section 2 we give proofs of the above results, and in Section 3 we make some concluding remarks.

\section{Proofs of the Main Results}

We denote by $N(v),$ $N[v]$ and ${\rm deg}(v)$ the neighborhood, closed neighborhood and degree of a vertex $v$ respectively.
By definition, $N[v]=\{v\}\cup N(v).$ Given a graph $G,$ $\Delta(G)$
and $\overline{G}$ denote the maximum degree of $G$ and the complement of $G$ respectively.
For $S\subseteq V(G),$ we denote by $G[S]$ the subgraph of $G$ induced by $S.$
For two graphs $G$ and $H,$ $G\vee H$ denotes the {\it join} of $G$ and $H,$ which is obtained from the disjoint union $G+H$
by adding edges joining every vertex of $G$ to every vertex of $H.$ $P_n$ denotes the path of order $n$ (and hence of length
$n-1.$) As usual, $q K_2$ denotes the graph consisting of $q$ pairwise vertex-disjoint edges.

{\bf Notation.} For an even positive integer $n,$ the notation $K_n-PM$ denotes the graph obtained from the complete graph $K_n$
by deleting all the edges in a perfect matching of $K_n$; i.e., it is the complement of $(n/2)K_2.$

The following lemma is well-known [2,  pp.12-13], but its hamiltonian part is usually not stated.

{\bf Lemma 6.} {\it Let $k$ and $n$ be integers with $1\le k\le n-1.$ Then there exists a $k$-regular graph of order $n$ if and
only if $kn$ is even. If $kn$ is even and $k\ge 2,$ then there exists a hamiltonian $k$-regular graph of order $n.$ }

{\bf Lemma 7.} {\it Let $d$ and $n$ be integers with $1\le d\le n-1,$ and
let $f(n,d)$ be the maximum size of a graph of order $n$ with maximum degree $\le d.$ Then
$$
f(n,d)=\begin{cases}(nd-1)/2\quad if\,\,\,both\,\,\,n\,\,\,and\,\,\,d\,\,\,are\,\,\,odd;\\
                     nd/2 \quad otherwise.
       \end{cases}
$$ }

{\bf Proof.} If at least one of $n$ and $d$ are even, then by Lemma 6 there exists a $d$-regular graph of order $n.$ Hence the obvious
upper bound $nd/2$ can be attained.

If both $n$ and $d$ are odd, then such a graph has at most $n-1$ vertices with degree $d,$ since the number of odd vertices in any graph is even.
Hence $f(n,d)\le ((n-1)d+(d-1))/2=(nd-1)/2.$ Next we show that this upper bound can be attained. If $d=1,$ the graph $((n-1)/2)K_2+K_1$
attains the upper bound. Now suppose $d\ge 2.$ By Lemma 6, there exists a $d$-regular graph $G$ of order $n-1$ containing a matching
$M=\{x_iy_i|i=1,\ldots,(d-1)/2\}$ of size $(d-1)/2.$ In $G,$ delete all the edges in $M,$ add a new vertex $v$ and join the edges $vx_i,\, vy_i$
for $i=1,\ldots,(d-1)/2.$ Thus we obtain a graph with degree sequence $d,d,\ldots,d,d-1$ which has size $(nd-1)/2.$  $\Box$

{\bf Proof of Theorem 1.} A graph $G$ contains no $K_{1,\, p}$ if and only if $\Delta (G)\le p-1.$ By Lemma 7, ${\rm ex}(n,K_{1,\, p})=(n(p-1)-1)/2.$
Note that a graph with degree sequence $d_1\ge d_2\ge\ldots\ge d_n$ contains exactly one copy of $K_{1,\, p}$ if and only if $d_1=p$ and
$d_2\le p-1.$ Hence a graph of order $n$ and size ${\rm ex}(n,K_{1,\, p})+1$ contains exactly one copy of $K_{1,\, p}$ if and only if its
degree sequence is $p,p-1,\ldots,p-1.$

Our assumption in Theorem 1 implies $p-1\ge 3.$ By Lemma 6, there exists a hamiltonian $(p-1)$-regular graph $R$ of order $n-1.$
Obviously $R$ contains a matching $M=\{x_iy_i|i=1,\ldots,p/2\}$ of size $p/2.$ Deleting all the edges in $M,$ adding a new vertex $v,$ and joining
the edges $vx_i$ and $vy_i$ for $i=1,\ldots,p/2,$ we obtain a graph $H$ with degree sequence $p,p-1,\ldots,p-1.$ $H$ has order $n,$ size
${\rm ex}(n,K_{1,\, p})+1$ and contains exactly one copy of $K_{1,\, p}.$ $\Box$

{\bf Lemma 8.} {\it If $p$ is an even positive integer, then ${\rm ex}(p+2, B_p)=p(p+2)/2$ and ${\rm ex}(p+3, B_p)=p(p+4)/2.$
If $p$ is an odd positive integer, then ${\rm ex}(p+2, B_p)=(p+1)^2/2$ and ${\rm ex}(p+3, B_p)=(p+1)(p+3)/2.$ }

{\bf Proof.} The proofs for the four Tur\'{a}n numbers have the same pattern, but these results hold for different reasons.
Let $G$ be a graph of order $n$ and size $e$ with vertices $v_1,\ldots,v_n$ such that ${\rm deg}(v_i)=d_i,$ $i=1,\ldots,n$
and $d_1\ge\ldots\ge d_n.$ These notations will be used throughout the proof. We assign different values to the order $n$ and size
$e$ in different cases.

(1) ${\rm ex}(p+2, B_p)$ for even $p.$ Suppose $n=p+2$ and $e\ge (p(p+2)/2)+1.$ Then $\sum_{i=1}^n d_i=2e\ge p(p+2)+2.$
If $d_2\le p,$ then $\sum_{i=1}^n d_i\le (p+1)+(p+1)p=(p+1)^2< p(p+2)+2,$ a contradiction. Hence $d_2\ge p+1,$ implying that $d_1=d_2=p+1.$
Then $G$ has $p$ triangles sharing the common edge $v_1v_2.$ Thus $G$ contains $B_p.$ This shows that ${\rm ex}(p+2, B_p)\le p(p+2)/2.$

On the other hand, the graph $G_1=K_{p+2}-PM$ has order $p+2$ and size $p(p+2)/2,$ and does not contain $B_p.$ In fact, every edge of $G_1$
lies in exactly $p-2$ triangles. To see this, let $uv$ be an edge of $G_1$ and let $ux$ and $vy$ be the two edges in the perfect matching.
Then $uvzu$ is a triangle of $G_1$ if and only if $z\not\in \{u,v,x,y\},$ and consequently there are $(p+2)-4=p-2$ choices for $z.$
$G_1$ yields ${\rm ex}(p+2, B_p)\ge p(p+2)/2.$ Combining this inequality with the proved reverse inequality we obtain
${\rm ex}(p+2, B_p)=p(p+2)/2.$

(2) ${\rm ex}(p+3, B_p)$ for even $p.$ Suppose $n=p+3$ and $e\ge (p(p+4)/2)+1.$ Then $\sum_{i=1}^n d_i=2e\ge p(p+4)+2.$
We distinguish two cases.

Case 1. $d_1=p+2.$ If $d_2\le p,$ then $\sum_{i=1}^n d_i\le (p+2)+(p+2)p=(p+1)(p+2)< p(p+4)+2,$ a contradiction. Hence $d_2\ge p+1.$
Let $v_1,w_1,w_2,\ldots,w_p$ be $p+1$ distinct neighbors of $v_2.$ Then $G$ contains the $p$ triangles $v_1v_2w_iv_1,$ $i=1,\ldots,p$
sharing the common edge $v_1v_2.$ Thus $G$ contains $B_p.$

Case 2. $d_1\le p+1.$ If $d_{n-1}\le p,$ then $\sum_{i=1}^n d_i\le (n-2)(p+1)+2p=(p+1)^2+2p<p(p+4)+2,$ a contradiction. Hence
$d_1=d_2=\cdots=d_{n-1}=p+1.$ If $d_n\le p-1,$ then again we have $\sum_{i=1}^n d_i\le(p+2)(p+1)+p-1< p(p+4)+2,$ a contradiction.
Thus $d_n\ge p.$ But $d_n\not= p+1,$ since otherwise $G$ would be a $(p+1)$-regular graph of odd degree $p+1$ and odd order $p+3,$
which is impossible by Lemma 6. We must have $d_n=p.$ Without loss of generality, suppose $N(v_n)=\{v_1,\ldots,v_p\}.$ It follows that
$N(v_{n-2})=\{v_1,\ldots,v_{n-3},v_{n-1}\}$ and  $N(v_{n-1})=\{v_1,\ldots,v_{n-2}\}.$ Clearly $G$ has the $p$ triangles $v_{n-1}v_{n-2}v_iv_{n-1}$
for $i=1,2,\ldots,p$ sharing the common edge $v_{n-1}v_{n-2}.$ Hence $G$ contains $B_p.$ This shows that ${\rm ex}(p+3, B_p)\le p(p+4)/2.$

The graph $G_2=\overline{K}_3\vee (K_p-PM)$ has order $p+3$ and size $p(p+4)/2.$ Next we show that $G_2$ does not contain $B_p.$
Let $V_1$ be the vertex set of the subgraph of $G_2$ isomorphic to $\overline{K}_3$ and let $V_2$ be the vertex set of the subgraph of $G_2$
isomorphic to $K_p-PM.$ Let $xy$ be an edge of $G_2.$ If $x,y\in V_2,$ then $xy$ lies in exactly $(p-4)+3=p-1$ triangles; if $x\in V_1$ and
$y\in V_2,$ then $xy$ lies in exactly $p-2$ triangles. $G_2$ yields ${\rm ex}(p+3, B_p)\ge p(p+4)/2,$ which, combined with the proved reverse
inequality, shows that ${\rm ex}(p+3, B_p)=p(p+4)/2.$

(3) ${\rm ex}(p+2, B_p)$ for odd $p.$ Suppose $n=p+2$ and $e\ge ((p+1)^2/2)+1.$ Then $\sum_{i=1}^n d_i=2e\ge (p+1)^2+2.$
If $d_2\le p,$ then $\sum_{i=1}^n d_i\le (p+1)+(p+1)p=(p+1)^2<(p+1)^2+2,$ a contradiction. Hence $d_1=d_2=p+1.$ Then $G$ contains $p$ triangles
sharing the common edge $v_1v_2,$ which form a $B_p.$ This shows that ${\rm ex}(p+2, B_p)\le (p+1)^2/2.$

The graph $G_3=K_1\vee (K_{p+1}-PM)$ has order $p+2$ and size $(p+1)^2/2.$ Let $f$ be an edge of $G_3.$ If $f$ is incident with the vertex of $K_1,$
then $f$ lies in exactly $p-1$ triangles; if $f$ is an edge of $K_{p+1}-PM,$ then $f$ lies in exactly $p-2$ triangles. Thus $G_3$ contains no $B_p.$
$G_3$ shows that ${\rm ex}(p+2, B_p)\ge (p+1)^2/2.$ This inequality, together with the reverse inequality, proves that
${\rm ex}(p+2, B_p)=(p+1)^2/2.$

(4) ${\rm ex}(p+3, B_p)$ for odd $p.$ Finally suppose $n=p+3$ and $e\ge ((p+1)(p+3)/2)+1.$ Then $\sum_{i=1}^n d_i=2e\ge (p+1)(p+3)+2.$
If $d_2\le p+1,$ then $\sum_{i=1}^n d_i\le (p+2)+(p+2)(p+1)=(p+2)^2<(p+1)(p+3)+2,$ a contradiction. Hence $d_1=d_2=p+2.$
Then the edge $v_1v_2$ lies in $p+1$ triangles, implying that $G$ contains a $B_{p+1}$ and hence a $B_p.$ This shows that
${\rm ex}(p+3, B_p)\le (p+1)(p+3)/2.$

On the other hand, the graph $G_4=K_{p+3}-PM$ has order $p+3$ and size $(p+1)(p+3)/2.$ Every edge of $G_4$ lies in exactly $p-1$ triangles.
Hence $G_4$ contains no $B_p.$ $G_4$ yields ${\rm ex}(p+3, B_p)\ge (p+1)(p+3)/2.$ We thus conclude that ${\rm ex}(p+3, B_p)=(p+1)(p+3)/2.$ $\Box$

{\bf Proof of Theorem 2.} By Lemma 8, ${\rm ex}(p+2, B_p)=p(p+2)/2.$ Let $G_5=K_2\vee (K_p-PM);$ i.e., $G_5$ is the complement of the graph $(p/2)K_2+\overline{K}_2.$ Note that $G_5$ is obtained from the graph $G_1$ in the proof
of Lemma 8 by adding one edge: $G_5=G_1+f.$ Then $G_5$ has order $p+2$ and size $(p(p+2)/2)+1.$ $f$ is the unique edge of $G_5$ which lies in
$p$ triangles. Hence $G_5$ contains exactly one copy of $B_p.$

Conversely, let $Y$ be a graph of order $p+2$ and size $(p(p+2)/2)+1$ which contains exactly one copy of $B_p.$ Let $uv$ be the unique edge
of $Y$ that lies in exactly $p$ triangles. Then ${\rm deg}(u)={\rm deg}(v)=p+1.$ Since $Y$ contains only one copy of $B_p,$ every vertex of $Y$
other than $u$ and $v$ has degree at most $p.$ The degree sum of $Y$ is $p(p+2)+2,$ implying that the degree sequence of $Y$ must be
$p+1,p+1,p,p,\ldots,p.$ It follows that $Y$ is the complement of $(p/2)K_2+\overline{K}_2.$ Hence $Y=G_5.$

Next we consider the order $p+3.$ By Lemma 8, ${\rm ex}(p+3, B_p)=p(p+4)/2.$ Let $G_6$ be the complement of $(p/2)K_2+P_3$ where $P_3$ is the path
of order $3.$ Note that $G_6$ is obtained from the graph $G_2$ in the proof of Lemma 8 by adding one edge to the subgraph $\overline{K}_3$:
$G_6=(K_1+K_2)\vee (K_p-PM).$ Then $G_6$ has order $p+3$ and size $(p(p+4)/2)+1.$ Let $h$ be the edge of $G_6$ corresponding to $K_2;$ i.e.,
the edge added to $G_2.$ It is easy to check that $h$ is the unique edge of $G_6$ that lies in $p$ triangles. Hence $G_6$ contains exactly
one copy of $B_p.$

Conversely, let $Z$ be a graph of order $p+3$ and size $(p(p+4)/2)+1$ which contains exactly one copy of $B_p.$
Let $v_1,\ldots,v_{p+3}$ be the vertices of $Z$ such that $v_1$ and $v_2$ are adjacent, $N(v_1)\cap N(v_2)=\{v_3,v_4,\ldots,v_{p+2}\}$ and $v_{p+3}\not\in N(v_1)\cap N(v_2).$ Thus we assumed that $v_1v_2$ is the unique edge of $Z$ that lies in $p$ triangles.
If ${\rm deg}(v_i)=p+2$ for some $i$ with $3\le i\le p+2,$ then $Z$ would contain at least three copies of $B_p,$ a contradiction.
Hence ${\rm deg}(v_i)\le p+1$ for $3\le i\le p+2.$ We also have ${\rm deg}(v_{p+3})\le p,$ since otherwise $Z$ would contain at least two copies of
$B_p.$

We assert that $v_{p+3}\not\in N(v_1)\cup N(v_2).$ To the contrary, without loss of generality
suppose $v_{p+3}$ is adjacent to $v_1.$ First consider the case $p=2.$ Then $v_{p+3}=v_5$ cannot be adjacent to any of $v_3$ and $v_4,$
and also $v_3$ and $v_4$ are not adjacent, since otherwise $Z$ would contain at least two copies of $B_2.$ But then $Z$ has size
$6<7=(2(2+4)/2)+1,$ a contradiction. Next assume $p\ge 4.$ If ${\rm deg}(v_{p+3})\le p-2,$ then the degree sum of $Z$ is at most
$(p+2)+(p+1)(p+1)+(p-2)=p^2+4p+1<p(p+4)+2,$ the degree sum of $Z,$ which is a contradiction. Thus ${\rm deg}(v_{p+3})\ge p-1,$ implying that
$v_{p+3}$ has at least $p-2$ neighbors among the vertices $v_3,\ldots,v_{p+2.}$ At least one, say $v_j,$ of these $p-2$ neighbors of $v_{p+3}$
has degree $\ge p+1,$ since otherwise the degree sum of $Z$ is at most
$$
(p+2)+3(p+1)+(p-2)p+p=p^2+3p+5<p(p+4)+2,
$$
a contradiction. But now the edge $v_1v_j$ lies in $p$ triangles, yielding
another copy of $B_p,$ a contradiction. This proves that $v_{p+3}\not\in N(v_1)\cup N(v_2).$ It follows that ${\rm deg}(v_1)={\rm deg}(v_2)=p+1.$

Summarizing the above analysis we obtain ${\rm deg}(v_1)={\rm deg}(v_2)=p+1,$ ${\rm deg}(v_i)\le p+1$ for $3\le i\le p+2$ and
${\rm deg}(v_{p+3})\le p.$ These restrictions on the degrees, together with the condition that the degree sum of $Z$ is $p(p+4)+2,$
imply that the degree sequence of $Z$ must be $p+1,p+1,\ldots, p+1,p.$ It remains to show that $G_6$ is the only graph with such a degree sequence,
 implying that $Z=G_6.$ The complement of such a graph has degree sequence $2,1,\ldots,1.$ Clearly $(p/2)K_2+P_3$ is the only graph with degree sequence $2,1,\ldots,1.$ This completes the proof. $\Box$

{\bf Proof of Corollary 3.} Theorem 2 covers the case when $k$ is even, and Theorem 1 covers the case when $k$ is odd.  $\Box$

{\bf Proof of Theorem 4.} Let $G$ be a graph of order $n$ and size ${\rm ex}(n, B_p)+1$ with vertices $v_1,\ldots,v_n$ such that ${\rm deg}(v_i)=d_i,$ $i=1,\ldots,n$ and $d_1\ge\ldots\ge d_n.$

First we suppose $n=p+2.$ By Lemma 8, ${\rm ex}(n, B_p)=(p+1)^2/2.$ Then the degree sum of $G$ is $\sum_{i=1}^n d_i=(p+1)^2+2.$ If $d_3\le p,$
then $\sum_{i=1}^n d_i\le 2(p+1)+p\times p <(p+1)^2+2,$ a contradiction. Thus we have $d_1=d_2=d_3=p+1.$ Each of the three edges $v_1v_2,$
$v_2v_3$ and $v_3v_1$ lies in $p$ triangles. Hence $G$ contains at least $3$ copies of $B_p.$ The number $3$ is attained by the graph
$K_3\vee (K_{p-1}-PM),$ which has order $p+2$ and size $((p+1)^2/2)+1,$ and contains exactly $3$ copies of $B_p.$

Next we suppose $n=p+3.$ By Lemma 8, ${\rm ex}(n, B_p)=(p+1)(p+3)/2.$ We have $\sum_{i=1}^n d_i=(p+1)(p+3)+2.$ If $d_2\le p+1,$ then
$\sum_{i=1}^n d_i\le (p+2)+(p+2)(p+1)=(p+2)^2< (p+1)(p+3)+2,$ a contradiction. Hence $d_1=d_2=p+2.$ We distinguish two cases.

Case 1. $d_3=p+2.$ Then each of the three edges $v_1v_2,$ $v_2v_3$ and $v_3v_1$ lies in $p+1$ triangles, implying that $G$ contains $3$ copies
of $B_{p+1}.$ Since one copy of $B_{p+1}$ contains $p+1$ copies of $B_p,$ $G$ contains at least $3(p+1)$ copies of $B_p.$

Case 2. $d_3\le p+1.$ The degree sum $(p+1)(p+3)+2$ of $G$ requires that $d_i=p+1$ for each $i=3,4,\ldots,n.$ Now in $G,$ the edge $v_1v_2$
lies in $p+1$ triangles, yielding $p+1$ copies of $B_p;$ for every $i=3,\ldots,p+3,$ each of the two edges $v_iv_1$ and $v_iv_2$ lies in
$p$ triangles, yielding $2(p+1)$ copies of $B_p.$ Altogether, $G$ contains at least $3(p+1)$ copies of $B_p.$

The number $3(p+1)$ is attained by the graph $K_2\vee (K_{p+1}-PM),$ which has order $p+3$ and size $((p+1)(p+3)/2)+1,$ and contains exactly
$3(p+1)$ copies of $B_p.$ This completes the proof. $\Box$

{\bf Proof of Theorem 5.} Clapham, Flockhart and Sheehan [3, p.36] determined the values ${\rm ex}(n,C_4)$ for all $n$ up to $21.$ The eight values we need are
${\rm ex}(6,C_4)=7,$ ${\rm ex}(7,C_4)=9,$ ${\rm ex}(8,C_4)=11,$ ${\rm ex}(9,C_4)=13,$ ${\rm ex}(10,C_4)=16,$ ${\rm ex}(11,C_4)=18,$
${\rm ex}(12,C_4)=21,$ and ${\rm ex}(13,C_4)=24.$

The graph in Figure 1(a) has order 6 and size $8={\rm ex}(6,C_4)+1$ and contains exactly one copy of $C_4;$ the graph in Figure 1(b) has order 7 and
size $10={\rm ex}(7,C_4)+1$ and contains exactly one copy of $C_4;$ the graph in Figure 1(c) has order 8 and size $12={\rm ex}(8,C_4)+1$ and contains exactly one copy of $C_4.$
\vskip 3mm
\par
 \centerline{\includegraphics[width=4.0in]{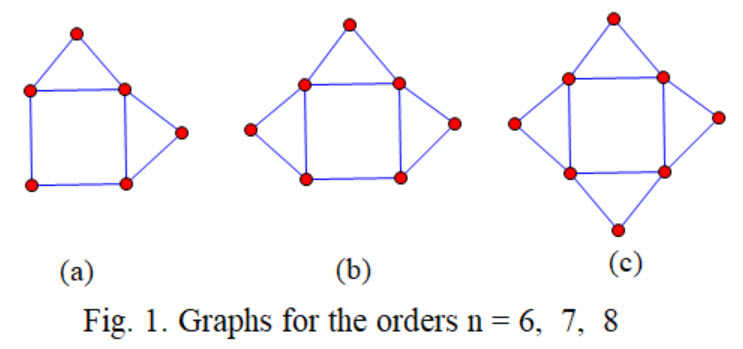}}
\par
The graph in Figure 2 has order 9 and size $14={\rm ex}(9,C_4)+1$ and contains exactly one copy of $C_4,$ the cycle $1,3,2,9,1.$
\vskip 3mm
\par
 \centerline{\includegraphics[width=2.9in]{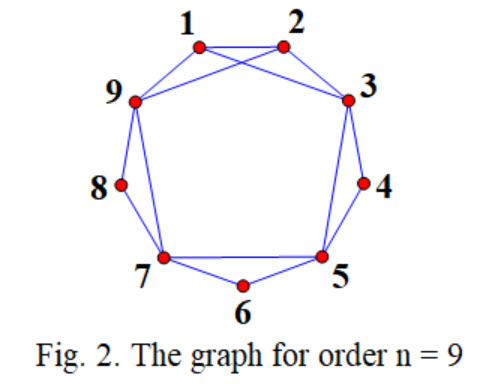}}
\par
The graph in Figure 3 has order 10 and size $17={\rm ex}(10,C_4)+1$ and contains exactly one copy of $C_4,$ the cycle $1,2,6,10,1.$
To find $4$-cycles in a graph of a small order, we check all the vertex pairs $i,j$ with $i<j$ and see whether the number
$|N(i)\cap N(j)|$ is at least $2.$ An inspection of Figure 3 gives $N(1)\cap N(6)=\{2,10\}$ and $N(2)\cap N(10)=\{1,6\},$
both  corresponding to the $4$-cycle $1,2,6,10,1.$
\vskip 3mm
\par
 \centerline{\includegraphics[width=2.9in]{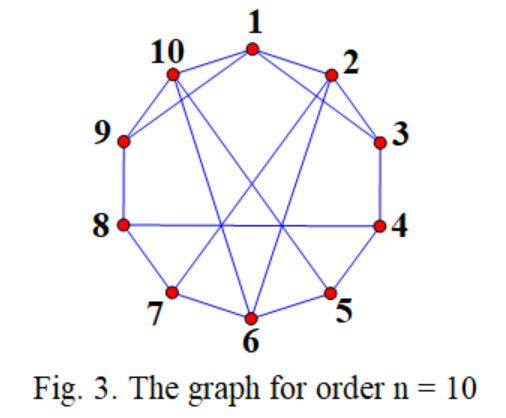}}
\par
The graph in Figure 4 has order 11 and size $19={\rm ex}(11,C_4)+1$ and contains exactly one copy of $C_4,$ the cycle $1,7,9,11,1.$
\vskip 3mm
\par
 \centerline{\includegraphics[width=2.9in]{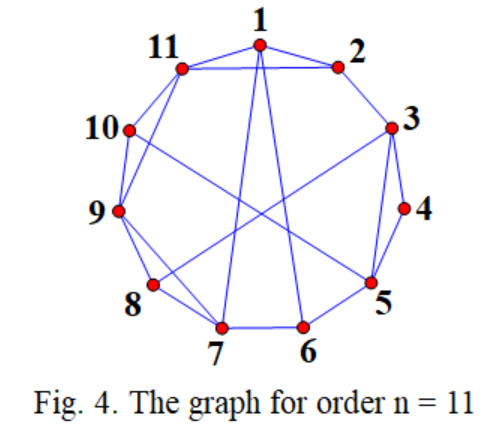}}
\par

Next we treat the case of order $n=12.$ Let $G$ be a graph of order $12$ and size $22={\rm ex}(12,C_4)+1.$ We first show that $G$ contains at least two copies of
$C_4.$

Let $v$ be a vertex of degree $\Delta (G).$  Denote $A=G[N(v)]$ and $B=G-N[v].$ We denote
by $[A,\,B]$ the set of edges of $G$ with one end in $A$ and the other end in $B.$

Let $\mathscr{S}$ be a set of graphs and let $n$ be a positive integer. The {\it Tur\'{a}n number of $\mathscr{S}$ for the order $n,$} denoted
${\rm ex}(n,\mathscr{S}),$ is the maximum size of a simple graph of order $n$ not containing any graph in $\mathscr{S}$ as a subgraph. The following
four facts can be easily verified.

{\bf Fact 1.}
$$
{\rm ex}(n,\{C_3,\,P_4,\,K_{1,3}\})=\begin{cases} 2k\quad if\,\,\,n=3k\,\,\,or\,\,\,n=3k+1;\\
                     2k+1 \quad if\,\,\,n=3k+2.
       \end{cases}
$$

{\bf Fact 2.} If $A$ contains any of the three graphs $C_3,\,P_4,\,K_{1,3},$ then $G[N[v]],$ and hence $G$ contains at least two
copies of $C_4.$

{\bf Fact 3.} Suppose that $B$ has order $k.$ If $|[A,\, B]| \ge k+2,$ then $G$ contains at least two copies of $C_4,$ and if
$|[A,\, B]| = k+1,$ then $G$ contains at least one copy of $C_4$ with exactly one vertex in $B.$

{\bf Fact 4.} (1) The minimum number of copies of $C_4$ in a graph of order $5$ and size $7$ is $2.$  (2) The minimum number of copies of $C_4$
in a graph of order $6$ and size $9$ is $3.$

We distinguish two cases according to $\Delta=\Delta (G).$ Since the average degree of $G$ is $11/3,$ we have $\Delta\ge 4.$ We use the notation
$e(H)$ to denote the size of a graph $H.$

Case 1. $5\le\Delta\le 11.$

Subcase 1.1. $7\le\Delta\le 11.$ We give only the proof of the case $\Delta =7,$ since the proofs for other values of $\Delta$ in this range
are the same. Now $B$ has order $4.$ If $e(B)=6,$ then it is $K_4$, which already contains three copies of $C_4.$ So suppose $e(B)\le 5.$
Also by Fact 3, we may suppose $|[A,\, B]| \le 5.$ Then $A$ has order $7$ and size at least $5.$ By Fact 1, $A$ contains
at least one of $C_3,\,P_4,\,K_{1,3},$ which implies that $G$ contains at least two copies of $C_4$ by Fact 2.

Subcase 1.2. $\Delta=6.$ $B$ has order $5.$ By Fact 4(1) we may suppose $e(B)\le 6,$ and by Fact 3 we may suppose $|[A,\, B]| \le 6.$
 If $|[A,\, B]| = 6,$ by Fact 3, $G$ has a copy of $C_4$
with exactly one vertex in $B.$ Now the subgraph $G[N[v]]$ has order $7$ and size at least $10;$ it contains at least one copy of $C_4$ since
${\rm ex}(7,C_4)=9.$ Thus $G$ contains at least two copies of $C_4.$ It remains to consider the case when $|[A,\, B]| \le 5.$
Then $A$ has order $6$ and size at least $5.$ By Fact 1, $A$ contains at least one of $C_3,$ $P_4$ and $K_{1,3}.$ Hence $G$ contains
at least two copies of $C_4$ by Fact 2.

Subcase 1.3. $\Delta=5.$ $B$ has order $6.$ By Fact 4(2) we may suppose $e(B)\le 8,$ and by Fact 3 we may suppose $|[A,\, B]| \le 7.$
First consider the case $e(B)=8.$ Then $B$ contains a copy of $C_4$ since ${\rm ex}(6,C_4)=7.$ If $|[A,\, B]| = 7,$ by Fact 3, $G$ contains
another copy of $C_4.$ If $|[A,\, B]| \le 6,$ then $G[N[v]]$ has order $6$ and size at least $8;$ it contains at least one copy of $C_4$ since ${\rm ex}(6,C_4)=7.$  Altogether $G$ contains at least two copies of $C_4.$ It remains to consider the case $e(B)\le 7.$ By distinguishing the
two cases when $|[A,\, B]| =7$ or $\le 6$ and using Facts 3 and 4, a similar analysis proves the result.

Case 2. $\Delta =4.$ Since ${\rm ex}(12,C_4)=21,$ $G$ contains a $4$-cycle $C.$ If $C$ has a vertex of degree $\le 3,$ deleting that vertex
we obtain a graph of order $11$ and size at least $19,$ which contains another $4$-cycle since ${\rm ex}(11,C_4)=18.$ Thus we may suppose that
each of the four vertices of $C$ has degree $4.$ Let $C=v,x,y,z,v$ and we continue using the above notations related to the vertex $v.$

Now $A$ has order $4,$ and $B$ has order $7.$ Since ${\rm ex}(7,C_4)=9,$ we may suppose that $e(B)\le 9;$ otherwise $B$ contains a $C_4$
which is clearly different from $C.$ By Fact 3, we may always suppose $|[A,\, B]| \le 8.$ By Facts 1 and 2, if $e(A)\ge 3,$ then $G$ contains
at least two copies of $C_4.$ Note that the above assumptions $e(B)\le 9$ and $|[A,\, B]| \le 8$ imply that $e(A)\ge 1.$
Next we suppose $1\le e(A)\le 2.$

Subcase 2.1. $v$ and $y$ are adjacent. First note that $A$ already has two edges $xy$ and $yz.$ If $|[A,\, B]|=8,$ by Fact 3,
$G$ contains a $C_4$ which is different from $C.$ It suffices to assume that $|[A,\, B]|\le 7.$ Since $e(G)=22,$ we must have $|[A,\, B]|=7$
and $e(B)=9.$ Then $y$ has exactly one neighbor in $B$ and every other vertex  of $A$ has exactly two neighbors in $B.$ Then for $G$ to avoid
having a second copy of $C_4$ under the condition $|[A,\, B]|=7,$ $B$ has size at most $8,$ a contradiction.

Subcase 2.2. $v$ and $y$ are nonadjacent and $e(A)=1.$ In this case we have $|[A,\, B]|=8$ and $e(B)=9.$
Assume $|[A,\, B]|=8.$ By analysing all the three possible places where the edge of $A$ lies, it can be checked that for $G$ to
avoid having a second copy of $C_4,$ $B$ has size at most $8,$ a contradiction.

Subcase 2.3. $v$ and $y$ are nonadjacent and $e(A)=2.$ Let $N(v)=\{x,z,p,q\}.$ For $G[N[v]]$ to avoid having a second copy of $C_4,$
there are essentially two possibilities for the two edges of $A:$ $E(A)=\{xz,\, pq\}$ or $E(A)=\{xp,\, zq\}.$ If $|[A,\, B]|=8$ then
$e(B)=8.$ But in both cases, for $G$ to avoid having a second copy of $C_4$ under the condition $|[A,\, B]|=8,$
 $B$ has size at most $7,$ a contradiction. If $|[A,\, B]|=7$ then $e(B)=9.$ But again in both cases, for $G$ to avoid having a second copy of $C_4$
 under the condition $|[A,\, B]|=7,$ $B$ has size at most $8,$ a contradiction.
We have proved that a graph of order $12$ and size $22$ has at least two $4$-cycles.

The graph in Figure 5 has order $12$ and size $22,$ and contains exactly two $4$-cycles: $3,4,5,9,3$ and $3,8,10,9,3.$ Thus the minimum number of
copies of $C_4$ in a graph of order $12$ and size $22$ is $2.$
\vskip 3mm
\par
 \centerline{\includegraphics[width=5.1in]{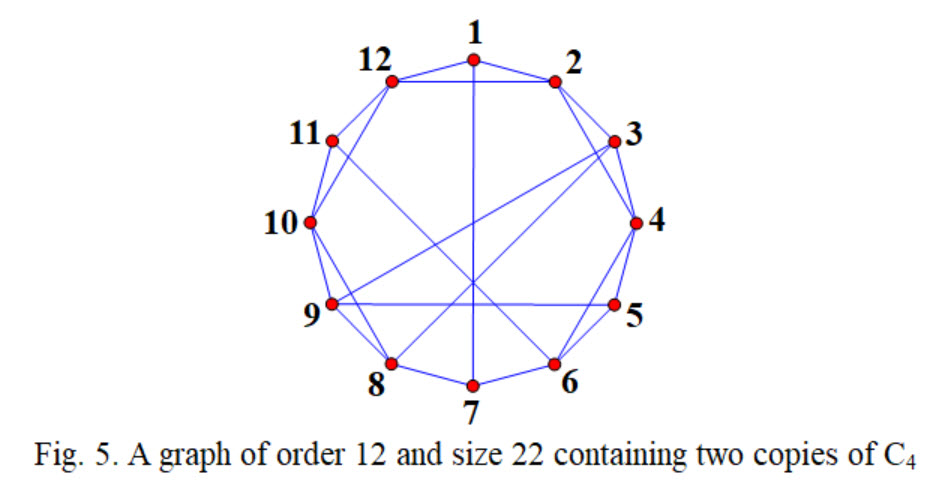}}
\par

Finally we consider the case of order $n=13.$ A graph of order $13$ and size $25={\rm ex}(13,C_4)+1$ has average degree $50/13<4.$
Hence it has a vertex of degree $\le 3.$ Deleting that vertex, we obtain a graph of order $12$ and size at least $22,$ which contains
at least two $4$-cycles by the proved result for the order $n=12.$ On the other hand, the graph in Figure 6 has order $13$ and size $25,$
and contains exactly two $4$-cycles: $1,2,10,12,1$ and $3,5,11,4,3.$
\vskip 3mm
\par
 \centerline{\includegraphics[width=5.3in]{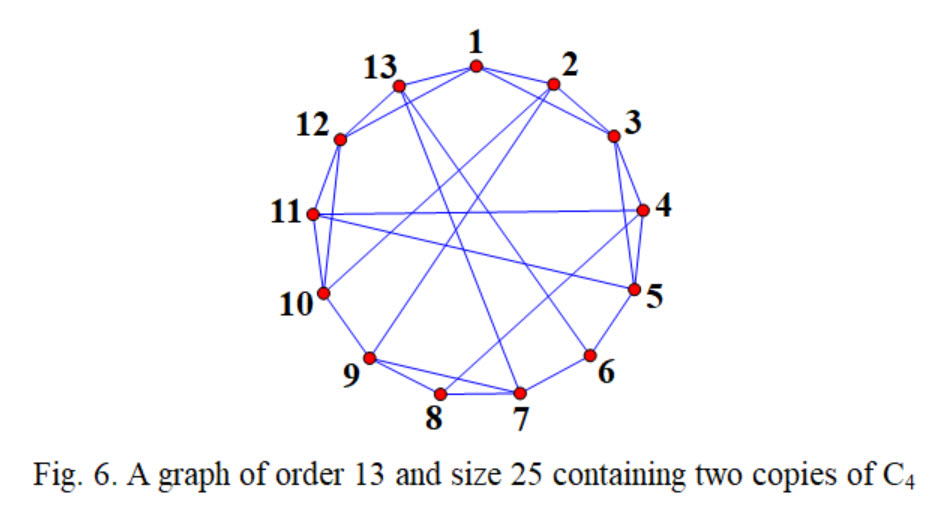}}
\par
Thus the minimum number of copies of $C_4$ in a graph of order $13$ and size $25$ is $2.$ This completes the proof. $\Box$

\section{Concluding Remarks}

We believe that Erd\H{o}s' intuition expressed in Problem 1 is almost true; i.e., for most graphs $H$ and most positive integers $n,$
every graph of order $n$ and size ${\rm ex}(n,H)+1$ contains at least two copies of $H.$ Now Problem 1 has been reduced to the following one.

{\bf Problem 2.} Determine all the pairs $(H,n),$ where $H$ is a graph and $n$ is a positive integer, such that there exists a graph
of order $n$ and size ${\rm ex}(n,H)+1$ which contains exactly one copy of $H.$

It is natural to ask whether Theorem 2 on the books $B_p$ can be extended to orders larger than $p+3.$ The answer is no in general.
A computer search gives the following information. (1) ${\rm ex}(6,B_2)=9,$ and the minimum number of copies of $B_2$ in a graph of
order $6$ and size $10$ is $2.$ (2) ${\rm ex}(7,B_2)=12,$ and the minimum number of copies of $B_2$ in a graph of order $7$ and size $13$ is $3.$
(3) ${\rm ex}(8,B_4)=21,$ and the minimum number of copies of $B_4$ in a graph of order $8$ and size $22$ is $6.$
(4) ${\rm ex}(9,B_4)=27,$ and the minimum number of copies of $B_4$ in a graph of order $9$ and size $28$ is $21.$

Finally we pose the following problem.

{\bf Problem 3.} Given positive integers $p$ and $n$ with $n\ge p+2,$ determine the Tur\'{a}n number ${\rm ex}(n,B_p).$

Lemma 8 solves the cases $n=p+2,\, p+3$ of Problem 3.

\vskip 5mm
{\bf Acknowledgements.} The authors are grateful to Professor Jie Ma from whose talk at ECNU based on the paper [9] they first learned
Problem 1. This research  was supported by the NSFC grants 11671148 and 11771148 and Science and Technology Commission of Shanghai Municipality (STCSM) grant 18dz2271000.

\end{document}